\input harvmac
\baselineskip=12pt plus 2pt minus 1pt
\noblackbox

\def\hom{{\rm Hom}}
\def\Hom{{\rm Hom}}
\def\GL{{\rm GL}}
\def\C{{\bf C}}
\def\Z{{\bf Z}}
\def\ra{\rightarrow}
\def\bfc{{\bf C}}
\def\bfz{{\bf Z}}
\def\bfr{{\bf R}}
\def\fuk{\cal F(\mir)}
\def\mir{\widetilde M}

\lref\syz{A. Strominger, S.-T. Yau, and E. Zaslow, ``Mirror
Symmetry is $T$-Duality,'' Nucl. Phys. {\bf B479} (1996) 243-259.}
\lref\vw{C. Vafa and E. Witten, ``On Orbifolds with Discrete Torsion,''
J. Geom. Phys. {\bf 15} (1995) 189.}
\lref\gh{P. Griffiths and J. Harris, {\sl Principles of Algebraic
Geometry,} Wiley \& Sons, New York, 1978.}
\lref\igusa{J. Igusa, {\sl Theta Functions,} Springer-Verlag, New York,
1972.}
\lref\helix{E. Zaslow, ``Solitons and Helices:  The Search for
a Math-Physics Bridge,'' Commun. Math. Phys. {\bf 175} (1996) 337-375.}
\lref\rfuk{K. Fukaya, Morse Homotopy, $A^\infty$-Category, and
Floer Homologies,'' in {\sl The Proceedings of the
1993 GARC Workshop on Geometry and Topology,} H. J. Kim, ed.,
Seoul National University;
``Floer Homology, $A^\infty$-Categories and
Topological Field Theory'' (notes by P. Seidel), Kyoto University
preprint Kyoto-Math 96-2.}
\lref\fmw{R. Friedman, J. Morgan, and E. Witten, ``Vector Bundles and
$F$ Theory,'' hep-th/9701162.}
\lref\kont{M. Kontsevich, ``Homological Algebra of Mirror Symmetry,''
Proceedings of the 1994 International Congress of Mathematicians {\bf I},
Birk\"auser, Z\"urich, 1995, p. 120;
alg-geom/9411018.}
\lref\rh{F. Reese Harvey, {\sl Spinors and Calibrations,} Academic Press,
New York, 1990.}
\lref\bbs{K. Becker, M. Becker, and A. Strominger, ``Fivebranes, Membranes,
and Non-Perturbative String Theory," hep-th/9507158.}
\lref\cdgp{P. Candelas, X. C. De La Ossa, P. Green, and L Parkes,
``A Pair of Calabi-Yau Manifolds as an Exactly Soluble Superconformal
Theory,'' Nucl.
Phys. {\bf B359} (1991) 21.}
\lref\gw{M. Gross and P. M. H. Wilson, ``Mirror Symmetry via 3-Tori for a Class
of Calabi-Yau Threefolds,'' alg-geom/9608004.}
\lref\mg{M. Gross, Special Lagrangian Fibrations I:  Topology," alg-geom/9710006.}
\lref\oog{D-Branes on Calabi-Yau Spaces and Their Mirrors," Nucl. Phys.
{\bf B477} (1996) 407-430.}
\lref\leung{N. C. Leung and C. Vafa, ``Branes and Toric Geometry,'' hep-th/9711013.}
\lref\atiyah{M. F. Atiyah, ``Vector Bundles over an Elliptic Curve,''
Proc. Lond. Math. Soc. {\bf VII} (1957) 414-452.}
\lref\hm{J. A. Harvey and G. Moore, ``On the Algebras of BPS States,''
hep-th/9609017.}
\lref\mor{D. R. Morrison, ``The Geometry Underlying Mirror Symmetry,''
alg-geom/9608006.}
\lref\gp{B. R. Greene and M. R. Plesser, ``Duality in Calabi-Yau
Moduli Space,'' Nucl. Phys. {\bf B338} (1990) 15-37.}
\lref\hl{R. Harvey and H. B. Lawson, Jr., ``Calibrated Geometries,''
Acta Math. {\bf 148} (1982) 47-157.}
\lref\hitch{N. Hitchin, ``The Moduli Space of Special Lagrangian
Submanifolds,'' dg-ga/9711002.}
\lref\fourmuk{C. Bartocci, U. Bruzzo, D. H. Ruip\'erez,
J. M. M. Porras, ``Mirror Symmetry on K3 Surfaces
via Fourier-Mukai Transform,'' alg-geom/9704023.}
\lref\kthree{D. Orlov, ``Equivalences of Derived Categories
and K3 Surfaces,'' alg-geom/9606006.}
\lref\rt{Y. Ruan and G. Tian, ``A Mathematical Theory of
Quantum Cohomology,'' J. Diff. Geom. {\bf 42} (1995) 259-367.}
\lref\wit{E. Witten, ``Mirror Manifolds and Topological Field
Theory,'' in {\sl Essays on Mirror Symmetry,} S.-T. Yau, ed.,
International Press, Hong Kong, 1992, pp. 120-159.}
\lref\ach{B.Acharya, ``A Mirror Pair of Calabi-Yau
Fourfolds in Type II String Theory,'' hep-th/9703029.}
\lref\stash{J. Stasheff, ``Homotopy Associativity of
$H$-Spaces.  I and II,'' Trans. Amer. Math. Soc {\bf 108}
(1963) 275-292 and 293-312.}
\lref\konttwo{M. Kontsevich, ``Deformation Quantization of
Poisson Manifolds, I,'' q-alg/9709040; S. Barannikov and M. Kontsevich,
``Frobenius Manifolds and Formality of Lie Algebras
of Polyvector Fields,'' alg-geom/9710032.}
\lref\mum{D. Mumford, {\sl Tata Lectures on Theta I,}
Birkh\"auser, Boston, 1983.}
\lref\ap{D. Arinkin, A. Polishchuk, ``Fukaya category and Fourier
transform,'' math.AG/9811023.}
\lref\pol{A. Polishchuk, ``Homological mirror symmetry with higher
products,'' math.AG/9901025.}

%-------------------
% title page
%-------------------
%
\pageno=275
\vsize=7.8in
\hsize=5in
\Title{}{\vbox{\hsize=5in\vskip 1in\centerline{\hskip -0.5in Categorical Mirror Symmetry:}
\vskip 0.2in
\centerline{\hskip -0.5in The Elliptic Curve}}}
\vskip -0.2in
\centerline{\hskip -0.5in { Alexander Polishchuk}\footnote{$^{\dagger}$}{email:
apolish@math.harvard.edu} {and
Eric Zaslow} \footnote{$^{\dagger\dagger}$}{email:
zaslow@math.harvard.edu}}
\vskip 0.1in
\centerline{\hskip -0.5in \it Department of Mathematics, Harvard University,
Cambridge, MA 02138, USA}
  
\vskip 0.3 in
\hsize=5in  
  
\centerline{\bf Abstract}
\vskip 0.1in
  
We describe an isomorphism of categories conjectured by Kontsevich.
If $M$ and $\widetilde{M}$ are mirror pairs then the conjectural equivalence
is between the derived category of coherent sheaves on $M$ and a
suitable version of Fukaya's category of Lagrangian submanifolds on
$\widetilde{M}.$  We prove this equivalence when $M$ is an elliptic
curve and $\widetilde{M}$ is its dual curve, exhibiting the dictionary
in detail.

\vskip 2in  
\rightline{previously published in ATMP 2 (1998) 443-470}
\Date{\rightline{\copyright 1999 International Press}}
\pageno=276
 
\newsec{Introduction and Summary}
  
Mirror symmetry -- equivalence between superconformal sigma models on
certain pairs of Calabi-Yau spaces -- has grown into an industry since its
discovery in 1989 \gp\
and the subsequent translation into enumerative
geometry in \cdgp.
Nevertheless, its origins remain mysterious.
In 1994, Kontsevich conjectured that mirror symmetry could be
interpreted as an equivalence of categories over mirror pairs of spaces
\kont.
While this would not explain mirror symmetry, it may be a way of
putting it on more algebraic grounds.  As it currently stands --
relating periods of Picard-Fuchs solutions
in an appropriate basis to Gromov-Witten invariants \rt\ --
mirror symmetry is rather difficult to define rigorously.
Further, it has been pointed out in \wit\ that the appropriate
moduli space for string theory should be more than just the K\"ahler and
complex moduli of the Calabi-Yau manifold.  Kontsevich introduced the
constructions considered here in part as a way of
naturally extending this moduli space (see also \konttwo).
 
\nopagenumbers
\headline={\ifodd\pageno\rightheadline \else\leftheadline\fi}
\def\rightheadline{\tenrm\hfil Alexander Polishchuk and Eric Zaslow  \hfil\folio}
\def\leftheadline{\tenrm\folio\hfil Categorical Mirror Symmetry:  The Elliptic Curve  \hfil}

All of the above work had been done without the knowledge of D-branes.
On a Calabi-Yau three-fold (in Type IIB), D-branes correspond to minimal
Lagrangian submanifolds with unitary local systems (by which we mean
flat $U(n)$ gauge bundles).  The conjecture of \syz\ states that these
objects are responsible for mirror symmetry.
For tori, K3, and certain orbifolds of tori, there is evidence
for such a description \vw, \ach, \gw, \leung.  We find
it appealing, then, that Kontsevich's
conjecture nicely incorporates these objects on the mirror side.
Indeed, the conjecture of \syz\ may now be reinterpreted
as a consequence of the correspondence between distinguished objects
in two equivalent categories.
  
Kontsevich's conjecture is that if $M$ and $\widetilde M$ are mirror pairs then
${\cal D}^b(M),$
the bounded derived category of coherent sheaves on $M,$
is equivalent to the
derived category of a suitable enlargement of the category ${\cal F}(\mir)$
of minimal Lagrangian submanifolds on $\widetilde M$
with unitary local systems:
${\cal D}^b(M) \cong {\cal D}^b({\cal F}(\mir)).$
${\cal F}(\mir)$ was introduced by Fukaya in \rfuk\ and refined by
Kontsevich \kont\ for the conjecture.
${\cal D}^b({\cal F}(\mir))$ denotes the derived category of the
$A^\infty$-category ${\cal F}(\mir),$ constructed in \kont\
using twisted complexes.  However, in the case of elliptic curves
one doesn't need such a construction.
Instead, we form a category
${\cal F}^0(\mir)$ from ${\cal F}(\mir)$ and prove\foot{In
the notation of \kont, ${\cal F}^0(\widetilde{M})$
would be called $H({\cal F}(\widetilde{M})).$}
$${\cal D}^b(M) \cong {\cal F}^0(\mir)$$
for the elliptic curve.
In general, it seems that the category ${\cal D}^b({\cal F}(\mir))$
has to be further localized to have a chance at equivalence with
${\cal D}^b(M).$
  
The torus is the first non-trivial check of this conjecture \kont.  Here
$M$ = $E_\tau$ is an elliptic curve with modular parameter $\tau$
and $\mir$ = $\widetilde E^\rho$ is a torus (the complex structure is irrelevant here)
with complexified K\"ahler parameter $\rho = iA + b,$ where $A$ is the
area and $b$ defines a class
in $H^2(\mir; {\bf R})/H^2(\mir;{\bf Z}).$
The mirror map is
$$\rho \leftrightarrow \tau$$
(due to the positivity of $A$ and the periodicity of $b,$ the more
natural parameter is perhaps $q = \exp(2\pi i \rho)).$
The equivalence that we detail in this simple case is non-trivial
and involves relations among theta functions and sections of
higher rank bundles, as well.  Note that our equivalence provides some
$A^\infty$-extension of the derived category of coherent sheaves
on an elliptic curve. This structure is discussed from the point of view
of complex geometry in \pol.
  
In the next section, we describe the
derived category of coherent sheaves
on a manifold, and then specifically on an elliptic curve.
Unfortunately, the elliptic curve is the only Calabi-Yau for which
${\cal D}^b$ is so well understood
(the case of K3 is discussed in \kthree).  In section three, we discuss
$\cal F,$ Kontsevich's generalization of Fukaya's category;
the category ${\cal F}^0;$ and
its definition for elliptic curves.
In section four we address the equivalence
\eqn\mainthm{{\cal D}^b(E_\tau) \cong {\cal F}^0(\widetilde E^\tau)}
informally, in the simplest case, giving
explicit examples along the way.  The general proof of \mainthm\
is our {\bf Main Theorem}, and comprises
section five.  Concluding remarks and speculations may be found
in section six.  A table of contents is listed below.
A physicist
interested in getting a gist for the equivalence may wish to
begin with the simple example of section four.
\bigskip 
%\vfill
%\eject
 
$$\eqalign{
&\hbox{\bf 1.  Introduction}\cr
&\hbox{\bf 2.  The Derived Category of Coherent Sheaves}\cr
&\qquad\hbox{\sl 2.1  Coherent Sheaves} \cr
&\qquad\hbox{\sl 2.2  The Derived Category}\cr
&\qquad\hbox{\sl 2.3  Vector Bundles on an Elliptic Curve}\cr
&\hbox{\bf 3.  Fukaya's Category, ${\cal F}(\mir)$}\cr
&\qquad\hbox{\sl 3.1  Definition}\cr
&\qquad\hbox{\sl 3.2  $A^\infty$ Structure}\cr
&\qquad\hbox{\sl 3.3  ${\cal F}^0(\mir)$}\cr
&\hbox{\bf 4.  The Simplest Example}\cr
&\hbox{\bf 5.  Proof of Categorical Equivalence}\cr
&\qquad\hbox{\sl 5.1  The Functor $\Phi$}\cr
&\qquad\hbox{\sl 5.2  $\Phi\circ m_2 = m_2 \circ \Phi$}\cr
&\qquad\hbox{\sl 5.3  Isogeny}\cr
&\qquad\hbox{\sl 5.4  The General Case}\cr
&\qquad\hbox{\sl 5.5  Extension to Torsion Sheaves}\cr
&\hbox{\bf 6.  Conclusions}
}$$
  
\newsec{The Derived Category of Coherent Sheaves}
  
\subsec{Coherent Sheaves}
  
Holomorphic vector bundles are locally free sheaves, by which we mean
that the space of sections
over a neighborhood
$U \ni z$ is ${\cal O}\oplus ...\oplus{\cal O} = {\cal O}^{\oplus r},$
where $r$ is the rank of the vector bundle.
Here ${\cal O}$ is the sheaf of holomorphic functions on $U.$
A coherent sheaf ${\cal S}$
admits a local presentation as an
exact sequence ${\cal O}^{\oplus p}\rightarrow
{\cal O}^{\oplus q}\rightarrow {\cal S}\rightarrow 0,$
and as a result a syzygy (the local version of projective resolution),
i.e. an exact sequence of sheaves
$${\cal E}_n\rightarrow{\cal E}_{n-1}\rightarrow ...\rightarrow{\cal E}_0
\rightarrow {\cal S}\rightarrow 0,$$
where the ${\cal E}_i$ are locally free.
The complex $({\cal E.})$ has homology $H_0({\cal E.}) = {\cal S},$
all others zero by exactness.  Therefore, we can think of ${\cal S}$
as roughly equivalent to a complex of locally free sheaves, if we
only take homology.  An example of a non-locally free sheaf is the
skyscraper sheaf ${\cal O}_{z_0}$
over a point on a one-dimensional manifold,
given by the exactness of the following (local) sequence:
\eqn\sky{{\cal O}\matrix{{}_{z-z_0}\cr\longrightarrow\cr{}}{\cal O}\longrightarrow
{\cal O}_{z_0}\longrightarrow 0.}
The map $(z-z_0)$ means multiplication by $z-z_0.$ Expanding functions
in a power series around $z_0,$ only the constants are not in the image
of this map.  Around every other point we can find a neighborhood such
that the map is surjective, and so we say that the {\sl stalk} over
$z_0$ is non-trivial, but trivial over every other point.
The rank
of the sheaf ${\cal O}_{z_0}$ is thus zero, but its degree is one,
since there is a global
section given by a choice of constant over $z_0.$
The kernel of
the map $z-z_0$ is the sheaf of holomorphic functions vanishing at
$z_0.$
More generally, replacing multiplication by $(z-z_0)$
above by multiplication by $(z-z_0)^n$, we get the definition of the
sheaf ${\cal O}_{nz_0}$ supported at the point $z_0$.
Every coherent sheaf on a complex curve is a direct sum of
vector bundles and these ``thickened" skyscraper sheaves.
The reason for this is that for any coherent sheaf,
we can form the
torsion part ${\cal F}_{tor},$ which is a direct sum of
thickened skyscrapers.  ${\cal F}_{tor}$ fits into the
exact sequence
$$0\longrightarrow {\cal F}_{tor}\longrightarrow {\cal F}
\longrightarrow{\cal G}\longrightarrow 0,$$
where ${\cal G} = {\cal F}/{\cal F}_{tor}$ is locally
free (a vector bundle).
However, one finds no nontrivial extensions by ${\cal F}_{tor},$
and so ${\cal  F} = {\cal G}\oplus {\cal F}_{tor}$ as claimed.
 
\subsec{The Derived Category}
 
The (bounded) derived category ${\cal D}^b(M)$ of coherent sheaves  on
$M$ is obtained from the
category of (bounded) complexes of coherent sheaves by inverting
quasi-isomorphisms -- i.e., morphisms inducing isomorphisms on
all  cohomology  sheaves of the complex.
In the case when $M$ is a complex
projective curve every object of ${\cal D}^b(M)$ is isomorphic
to the direct sum of objects of the form $F[n]$ where
$F$ is a coherent sheaf on $M$, $F[n]$ denotes  the  complex
with the only non-zero term (equal to $F$) in degree $-n$.
Thus, every object of ${\cal D}^b(M)$ is a direct sum
of objects of the form $F[n]$ where $F$ is either a vector
bundle or has support at a point.
 
\subsec{Vector Bundles on an Elliptic Curve}
 
We collect here some facts about theta functions, concentrating
on the elliptic curve \gh.
Line bundles on an elliptic curve $E = {\bf C}/\Lambda$ have a
simple description.  The lattice $\Lambda$ has $2$ generators,
one of which (by a rescaling) can be chosen to be $z \rightarrow z + 1.$
Then $X = {\bf C}^*/{\bf Z},$ where the coordinate on ${\bf C}^*$
is $u = \exp(2\pi i z),$ and the action of the lattice is $u\rightarrow qu,$
where $q = \exp(2\pi i \tau).$  Call $\pi'$ the projection
map $\pi': {\bf C}^* \rightarrow E.$  Now using the fact that
$H^1({\bf C}^*,{\cal O}) = H^2({\bf C}^*,{\cal O}) = 0,$ the long
exact sequence associated to the exact sheaf
sequence $0\rightarrow
{\bf Z}\rightarrow{\cal O}\rightarrow {\cal O}^*\rightarrow 0$
tells us that $H^1({\bf C}^*,{\cal O}^*)\cong H^2({\bf C}^*,{\bf Z}).$
This map is by definition the first Chern class, and therefore line bundles on ${\bf C}^*$
are determined by the first Chern class.  But since there are no
non-trivial two-forms on ${\bf C}^*,$ we learn that the pull-back of
any line bundle $L$ over $E$ is trivial over ${\bf C}^*.$
  
More generally, the pull-back of every vector bundle on $E_q$
to $\C^*$ is trivial since all vector bundles on $E$ are obtained from line
bundles by successive extensions. Thus, all vector bundles on
$E$ are obtained by the following construction.
Let $V$ be an $r$-dimensional
vector space and $A: \bfc^*\rightarrow GL(V)$
a holomorphic, vector-valued function.  We define the rank $r$
holomorphic vector bundle $F_q(V,A)$ on $E$ by taking
the quotient
$$F_q(V,A) = \bfc^*\times V/(u,v)\sim(uq,A(u)\cdot v).$$
The theory of vector bundles on elliptic curves can be rephrased as
a classification of such holomorphic functions
up to equivalence of $A(u)$ and $B(qu)A(u)B(u)^{-1}$
where $B:\C^*\ra\GL(V)$.
 
When $V=\bfc$ and $A = \varphi$ a holomorphic function, we call
$L_q(\varphi)$ the line bundle constructed in this way (or
simply $L(\varphi)$ if there is no ambiguity).
We define $L \equiv L_q(\varphi_0),$ where $\varphi_0(u) =
\exp(-\pi i \tau - 2\pi i z) = q^{-1/2}u^{-1}.$
The line bundle $L:=L(\phi_0)$ will play a distinguished
role in our considerations since the classical theta function
is a section of $L$. Every holomorphic line bundle on $E_q$
has form $t_x^*L\otimes L^{n-1}$ for some $n\in\Z$ and $x\in E_q$,
where $t_x$ is the map of translation by $x$ on $E_q$.
For a torus, the theta function has as indices $\tau,$
the modular parameter, and $(c',c'') \in {\bf R}^2/{\bf Z}^2,$
a translation parametrizing different line bundles of the
same degree.  The theta function is defined to be
$$\theta[c',c''](\tau,z) =
\sum_{m\in \bfz}\exp\{2\pi i[\tau(m+c')^2/2 + (m+c')(z+c'')]\}.$$
If $(c',c'')=(0,0),$ we will simply use the notation $\theta(\tau,z).$
Then the $n$ functions $\theta[a/n,0](n\tau,nz),$
$a\in \bfz/n\bfz,$ are the global sections
of $L^n.$
 
Now consider the natural $r$-fold
covering $\pi_r:E_{q^r}\rightarrow E_q$
which sends $u$ to $u.$  Clearly,
$\pi_r^{-1}(u) = \{ u,uq,uq^2,...,uq^{r-1}\}.$
We have the natural functors
of pull-back and push-forward associated with $\pi_r$.
More concretely, $\pi_r^*F_q(V,A)=F_{q^r}(V,A^r)$,
while $\pi_{r*}F_q(V,A)=F_q(V\otimes\C^r,\pi_{r*}A)$,
where $\pi_{r*}A(v\otimes e_i)=v\otimes e_{i+1}$
for  $i=1,\ldots,r-1$,  $\pi_{r*}A(v\otimes  e_{r})=Av\otimes
e_1$. It is easy to see that the functor of pull-back
commutes with tensor product and duality. Also one
has natural isomorphisms
$$\pi_{r*}(F_1\otimes\pi_r^*F_2)\cong\pi_{r*}(F_1)\otimes
F_2,$$
$$(\pi_{r*}(F))^*\cong\pi_{r*}(F^*),$$ and
$$H^0(E_q,\pi_{r*}(F))\cong H^0(E_{q^r},F).$$
It follows that
\eqn\homeq{\eqalign{\Hom(F_1,\pi_{r*}F_2)&\cong\Hom(\pi_r^*F_1,F_2),\cr
\Hom(\pi_{r*}F_1,F_2)&\cong\Hom(F_1,\pi_r^*F_2).}}
 
The following results will be useful in the sequel.
 
{\bf Proposition 1}
 
Every indecomposable bundle on $E_q$ is
isomorphic to a bundle of the form
$\pi_{r*}(L_{q^r}(\varphi)\otimes F_{q^r}(\C^k,\exp N)),$
where $N$ is a constant indecomposable nilpotent matrix,
$\varphi=t_x^*\varphi_0\cdot\varphi_0^{n-1}$
for some $n\in \bfz$ and $x \in \C^*$,
and $t_x$ represents translation by $x$.
 
This proposition can be deduced easily from the classification
of holomorphic bundles on elliptic curves due to M. Atiyah \atiyah.
 
{\bf Proposition 2}
 
Let $\varphi = t_x^*\varphi_0\cdot \varphi_0^{n-1},$ with $n>0.$
Then for any nilpotent endomorphism $N\in {\rm End}(V),$
there is a canonical isomorphism
$${\cal V}_{\varphi,N}: H^0\big(L(\varphi)\big)\otimes V \rightarrow H^0
\big(L(\varphi)\otimes
F(V,\exp N)\big).$$
 
Proof:  Put ${\cal V}_{\varphi,N}(f\otimes v) = \exp(DN/n) f\cdot v,$
where $D = -u{d\over du} = - {1\over 2\pi i}{d\over dz}.$
Now since $f$ is a section of $L=L(\varphi)$ we have
$f(qu) = \varphi(u) f(u).$
Now a section $w$ of $L\otimes F(V,N)$ is equivalent
to a holomorphic function of $u$ with values in $V$ obeying
$w(qu) = \varphi(u)\exp(N)\cdot w(u).$
First note that
$D$ commutes with the operation $u\rightarrow uq,$ so
$[Df](uq) = D(f(uq)).$ We also have
$D(\varphi f)= \varphi\cdot(D+n)f,$
which is easily checked.  As a result, we find
$$\eqalign{\exp(DN/n)f(qu)v &= \exp(DN/n)[\varphi(u) f(u) v]\cr
&= \sum_{k=0}^{\infty} {1\over k! n^k} D^k(\varphi f) N^k v \cr
&= \sum_{k=0}^{\infty} {1\over k! n^k} \varphi\cdot (D+n)^k f N^k v \cr
&= \varphi\cdot \exp[(D+n)N/n] fv \cr
&= \varphi\cdot \exp(N)[\exp(DN)\cdot fv],}$$
and the proposition is shown.
 
{\bf Proposition 3}
 
Let $\varphi_1 = t_{x_1}^*\varphi_0\cdot\varphi_0^{n_1-1},
\varphi_2 = t_{x_2}^*\varphi_0\cdot\varphi_0^{n_2-1},$
and let $N_i \in {\rm End}(V_i),$ $i = 1,2,$ be nilpotent endomorphisms.
Then
$$\eqalign{{\cal V}&_{\varphi_1,N_1}(f_1\otimes v_1)\circ
{\cal V}_{\varphi_2,N_2}(f_2\otimes v_2) = \cr
&{\cal V}_{\varphi_1\varphi_2,N_1 + N_2}
\left[\exp\left({n_2N_1- n_1 N_2\over n_1 + n_2}
{D\over n_1}\right) (f_1)
\exp\left({n_1 N_2 - n_2 N_1\over n_1 + n_2}
{D\over n_2}\right) (f_2) (v_1\otimes v_2)\right],}$$
where $N_1, N_2$ denote $N_1 \otimes 1$ and $1\otimes N_2$
respectively, on the right hand side, and
$\circ$ denotes the natural composition of
sections
$$\eqalign{
H^0\big( L(\varphi_1) \otimes F(V_1,\exp N_1)\big)& \otimes
H^0\big( L(\varphi_2) \otimes F(V_2,\exp N_2)\big) \rightarrow \cr
&H^0\big( L(\varphi_1\varphi_2)\otimes
F(V_1\otimes V_2,\exp(N_1\otimes 1 + 1\otimes N_2))\big).}$$
Recalling that $\exp {d\over dz}$ is the generator of translations,
we may write formally $\exp(N\cdot {d\over dz})f(z) = f(z+N).$
In this notation, the above formula looks like
$$\eqalign{{\cal V}(f_1\otimes v_1)\circ {\cal V}(f_2\otimes v_2) &=
{\cal V}\left(f_1\big(z + {n_1N_2 - n_2N_1 \over 2\pi i n_1(n_1+n_2)}\big)
f_2\big(z+{n_2N_1 - n_1N_2\over 2\pi i n_2 (n_1+n_2)}\big)(v_1\otimes v_2)
\right)\cr
&= {\cal V}\left( f_1\big(ue^{n_1N_2 - n_2N_1 \over n_1(n_1+n_2)}\big)
f_2\big( ue^{n_2N_1 - n_1N_2\over n_2 (n_1+n_2)}\big)(v_1\otimes
v_2)\right).}$$
 
The proof is straightforward.
 
\newsec{Fukaya's Category, $\cal F(\mir)$}
  
\subsec{Definition}
  
Let $\mir$ be a Calabi-Yau manifold with
its unique, Ricci-flat K\"ahler metric and K\"ahler
form $k$ (not just
the cohomology class, but the differential form).
Our main example will be when $\mir$ is a torus, $\widetilde{E},$
with flat metric $ds^2 = A(dx^2 + dy^2),$ so that its volume
is equal to $A.$  We also specify an element
$b$ in $H^2(\widetilde{E};{\bf R})/H^2(\widetilde{E};{\bf Z}).$
We define a complexified K\"ahler form $\omega \equiv b + ik.$
  
A category is given by a set of objects and composable morphisms between
objects.  $A^\infty$ categories have additional structures on the
morphisms, which we will discuss later in this section.
  
{\bf Objects:}  The objects of $\fuk$ are special Lagrangian submanifolds
of $\mir$ --
i.e. minimal Lagrangian submanifolds -- endowed with
flat bundles with monodromies having eigenvalues of unit
modulus\foot{Kontsevich only considered unitary local systems,
or flat $U(n)$ bundles.  We prove an equivalence of categories
involving a larger class of objects.  The Jordan blocks will
be related to non-stable vector bundles over the torus.  See
section five for details.}, and one additional structure
we will discuss momentarily.  We recall that a Lagrangian
submanifold is one on which $k$ restricts to zero.
Though Fukaya and Kontsevich
have taken the submanifolds to simply be Lagrangian, we will need the minimality condition as
well.\foot{A Lagrangian minimal (or ``special Lagrangian'') submanifold $L$
is a kind of calibrated submanifold \hl\ associated to the Calabi-Yau
form, $\Omega.$  This means that for a suitable complex
phase of $\Omega,$ we have ${\rm Re}(\Omega)\vert_L = {\rm Vol}_L,$ where
the volume form is determined by the induced metric.  Equivalently,
one has $k\vert_L = 0$ and ${\rm Im}(\Omega)\vert_L = 0.$}
This, plus Kontsevich's addition of the local system,
is also what one expects based on relations with
D-branes \bbs.  Thus, an object ${\cal U}_i$
is a pair:
$${\cal U}_i = ({\cal L}_i,{\cal E}_i).$$
  
In our example of a torus, the minimal submanifolds are just
minimal lines, or geodesics
(the Lagrangian property is trivially true for one-dimensional
submanifolds).  To define a closed submanifold in
${\bf R}^2/({\bf Z}\oplus{\bf Z})$ the slope of the line must be rational, so
can be given by a pair of integers $(p,q).$
There is another real datum needed, which is the point of interception
with the line $x=0$ (or $y=0$ if $p=0$).
In the easiest case, the rank of the unitary system is one, so
that we can specify a flat line bundle on the circle by simply
specifying the monodromy around the circle, i.e. a complex
phase $\exp(2\pi i \beta),$ $\beta \in {\bf R}/{\bf Z}.$
For a general local system of rank $r$ we can take $(p,q)$ to have
greatest common divisor equal to $r.$
  
The additional structure we need is the following.  A Lagrangian
submanifold $L$ of real dimension $n$ in a complex $n$-fold, $\mir,$
defines not only a map from $L$ to $\mir$ but also the Gauss map
from $L$ to $V,$ where $V$ fibers over $\mir$ with fiber at $x$
equal to the space of Lagrangian planes at $T_x\mir.$  The space of
Lagrangian planes has fundamental group equal to ${\bf Z},$ and we
take as objects special Lagrangian submanifolds
together with lfts of the Gauss map into
the fiber bundle over $\mir$
with fiber equal to the universal cover of the space of Lagrangian
planes.
  
For our objects, we thus require more than the slope, which
can be thought of as a complex phase with rational tangency,
and therefore as $$\exp{i \pi \alpha}.$$  We need a choice of
$\alpha$ itself.  Clearly, the ${\bf Z}$-degeneracy
represents the
deck transformations of the universal cover of the space of
slopes.  Shifts by integers correspond to shifts by grading
of the bounded complexes in the derived category.  There is no
natural choice of zero in ${\bf Z}.$
  
{\bf Morphisms:}  The morphisms $\hom({\cal U}_i,{\cal U}_j)$
are defined as
$$\hom({\cal U}_i,{\cal U}_j) =
{\bf C}^{\#\{{\cal L}_i\cap{\cal L}_j\}}\otimes
\hom({\cal E}_i,{\cal E}_j),$$
where the second ``Hom'' in the above represents homomorphisms
of vector spaces underlying the local systems at the points
of intersection.
There is a ${\bf Z}$-grading on the morphisms.  If $p$
is a point in ${\cal L}_i\cap{\cal L}_j$ then it has a
Maslov index $\mu(p) \in {\bf Z}.$\foot{A discussion of
the Maslov index may be found in chapter four of the first
reference in \rfuk.}
  
For our example, let us consider $\hom({\cal U}_i,{\cal U}_j),$
where the unitary systems have rank one, and where the lines
${\cal L}_i$ and ${\cal L}_j$ go through the origin.  Then
$\tan \alpha_i = q/p$ and $\tan \alpha_j = s/r,$ with $(p,q)$
and $(r,s)$ both relatively prime pairs.  For simplicity, one
can think of the lines as the infinite set of
parallel lines on the universal
cover of the torus, ${\bf R}\oplus{\bf R}.$  It is then easy to see
that there are
$$\vert ps - qr\vert $$
non-equivalent points of intersection.
Since $\hom({\bf C},{\bf C})$ is one dimensional, the monodromy is
specified by a single complex phase $T_i$ at each point of
intersection.  For rank $n$ local systems, $T_i$ would be
represented by an $n\times n$ matrix.
  
The ${\bf Z}$-grading on $\hom({\cal U}_i,{\cal U}_j)$ is constant
for all points of intersection in our example (they are all related
by translation).  If $\alpha_i,$ $\alpha_j$ are the real numbers
representing the logarithms of the slopes, as above, then
for $p \in {\cal L}_i\cap{\cal L}_j$ the
grading is given by
$$\mu(p) = - [\alpha_j - \alpha_i],$$
where the brackets represent the greatest integer.  Note that
$- [x] - [-x] = 1,$
which the Maslov index must obey for a one-fold.
The Maslov index is non-symmetric.  For
$p\in {\cal L}_i\cap {\cal L}_j$ in an $n$-fold, $\mu(p)_{ij}
+ \mu(p)_{ji} = n$.  The asymmetry is reassuring,
as we know that $\hom(E_i,E_j)$ is not symmetric in the case of
bundles.  It is the
extra data of the lift of the Lagrangian plane which allowed us
to define the Maslov index in this way.
  
Generally speaking, a category has composable morphisms which satisfy
associativity conditions.  This is not generally true for the category
$\cal F(\mir).$  However, we have instead on $\cal F(\mir)$ an additional
interesting structure making $\cal F(\mir)$ an $A^\infty$ category.
Associativity will hold cohomologically, in a way.  The equivalence
of categories that we will prove will involve a true category
${\cal F}^0(\mir),$ which we will construct
from $\cal F(\mir)$ in section {\sl 3.3.}
  
\subsec{$A^\infty$ Structure}
  
The category $\fuk$ has an $A^\infty$ structure, by which
we mean the
composable morphisms satisfy conditions analogous to those
of an $A^\infty$
algebra \stash.  An $A^\infty$ algebra is a generalization of a
differential, graded algebra.  Namely, it is a ${\bf Z}$-graded
vector space, with a degree one map, $m_1,$ which squares to zero
$((m_1)^2 = 0)$.  There are higher maps, $m_k: A^{\otimes k}\rightarrow A,$
as well.
  
\noindent{\bf Definition:}  An $A^\infty$ category, ${\cal F}$ consists of
  
$\bullet \;$ A class of objects ${\rm Ob}({\cal F});$
  
$\bullet \;$ For any two objects, $X, Y,$ a ${\bf Z}$-graded
abelian group of morphisms ${\rm Hom}(X,Y);$
  
$\bullet \;$ Composition maps
$$m_k: {\rm Hom}(X_1,X_2)\otimes{\rm Hom}(X_2,X_3)\otimes ...
{\rm Hom}(X_k,X_{k+1})\rightarrow {\rm Hom}(X_1,X_{k+1}),$$
$k\geq 1,$ of degree $2-k,$ satisfying the condition
$$\sum_{r=1}^{n}\sum_{s=1}^{n-r+1}(-1)^\varepsilon m_{n-r+1}\big(
a_1\otimes...\otimes a_{s-1}\otimes m_r (a_s \otimes ... \otimes
a_{s+r-1})\otimes a_{s+r}\otimes ... \otimes a_n\big) = 0$$
for all $n\geq 1,$ where $\varepsilon = (r+1)s + r(n + \sum_{j=1}^{s-1}{\rm deg}(a_j)).$
  
An $A^\infty$ category with one object is called an
$A^\infty$ algebra.  The first condition $(n=1)$
says that $m_1$ is a degree one
operator satisfying $(m_1)^2 = 0,$ so it is a co-boundary operator
which we can denote $d.$  The second
condition says that $m_2$ is a degree zero map
satisfying $d(m_2(a_1\otimes a_2)) = m_2(da_1 \otimes a_2) +
(-1)^{{\rm deg}(a_1)}m_2(a_1 \otimes da_2),$ so $m_2$ is a
morphism of complexes and induces a product on cohomologies.
The third condition says that $m_2$ is associative at the
level of cohomologies.
  
The $A^\infty$ structure on Fukaya's category is given by
summing over holomorphic maps (up to projective equivalence)
from the disc $D^2,$ which take the
components of the boundary $S^1 = \del D^2$ to the special Lagrangian
objects.  An element $u_j$ of ${\rm Hom}({\cal U}_j,{\cal U}_{j+1})$
is represented by a pair
$$u_j = t_j \cdot a_j,$$ where $a_j \in {\cal L}_j\cap {\cal L}_{j+1},$
and $t_j$ is a matrix in
${\rm Hom}({\cal E}_j\vert_{a_j},{\cal E}_{j+1}\vert_{a_j}).$
$$m_k(u_1\otimes ...\otimes u_k) = \sum_{a_{k+1} \in
{\cal L}_1\cap{\cal L}_{k+1}}
C(u_1,...,u_k,a_{k+1}) \cdot a_{k+1},$$
where (notation explained below)
\eqn\mcoef{C(u_1,...,u_k,a_{k+1}) = \sum_{\phi} \pm\,
{\rm e}^{2\pi i \int \phi^* \omega}\;
\cdot P {\rm e}^{\oint \phi^* \beta},}
is a matrix in
${\rm Hom}({\cal E}_1\vert_{a_{k+1}},{\cal E}_{k+1}\vert_{a_{k+1}}).$
Here we sum over (anti-)holomorphic maps $\phi: D^2\rightarrow
\mir,$ up to projective equivalence,
with the following conditions along the boundary:
there are $k+1$ points $p_j = {\rm e}^{2\pi i \alpha_j}$
such that $\phi(p_j) = a_j$ and
$\phi({\rm e}^{2\pi i \alpha})
\in {\cal L}_j$ for
$\alpha \in (\alpha_{j-1}, \alpha_{j}).$
In the above, $\omega = b + ik$ is the complexified K\"ahler form,
the sign is determined by an orientation in the space of
holomorphic maps
(it will always be positive for us),
and
$P$ represents a path-ordered integration, where $\beta$ is
the connection of the flat bundle along the local system
on the boundary.  Note that
in the case of all trivial local systems $(\beta \equiv 0),$ the weighting
is just the exponentiated complexified area of the map.
The path-ordered integral is defined by
$$P {\rm e}^{\oint \phi^* \beta}
= P {\rm e}^{\int_{\alpha_{k}}^{\alpha_{k+1}}\beta_k d\alpha}
\cdot t_k \cdot
P {\rm e}^{\int_{\alpha_{k-1}}^{\alpha^{k}}\beta_{k-1} d\alpha}
\cdot t_{k-1}\cdot ... \cdot
t_1 \cdot P {\rm e}^{\int_{\alpha_{k+1}}^{\alpha_{1}}\beta_1 d\alpha}
$$
(this formula is easily understood by reading right to left).
Following the integration along the boundary, we get a homomorphism
from ${\cal E}_1$ to ${\cal E}_{k+1}.$
The path ordering symbol is a bit superfluous above, since the
integral is one-dimensional, but we have retained it for exposition.
  
Fukaya has shown the $A^\infty$ structure of his category in \rfuk.  Our
modifications should not affect the proof.
 
\subsec{${\cal F}^0(\mir)$}
  
To define a true category which can be checked against the derived
category, we simply take $H^0$ of all the morphisms (recall that they
have the structure of complexes).
In our case of the elliptic curve, since $m_1 = d = 0$
the cohomology complex is the same as the original complex, and
so we simply take the degree zero piece of $\hom$ -- however, the
construction is valid in general.  Recall that $m_2$ was associative
at the level of cohomlogies.  This reduces to true associativity of
the composition of morphisms.
We call this category form in this way ${\cal F}^0(\mir).$
For our
example $m_2$
has degree zero, survives this restriction and is associative.
The higher $m$'s
are profected to zero in this category, so our equivalence will
be defined by constructing a dictionary of objects and checking
compatibility with $m_2.$
  
\newsec{The Simplest Example}
  
We now demonstrate the equivalence between the derived category of
coherent sheaves on an elliptic curve and the category
${\cal F}^0(\widetilde{E})$ on
the dual torus, in the simplest possible example.  A complete proof
will be given in the next section.
  
We begin with an elliptic curve
$E = {\bf C}/({\bf Z}\oplus \tau{\bf Z}).$
$\tau$ defines the complex structure of $E.$
Its mirror torus is $\widetilde{E} = {\bf R}^2/({\bf Z}\oplus{\bf Z}),$
with K\"ahler (metric) structure defined below, and related to
$\tau$ by the mirror map.
  
Let us first construct the dictionary.
Line bundles of degree (first Chern class) $d$
correspond to lines of slope $d.$ With this simple
rule, and the definition of the compositions
in both categories, we can try to compare the
map $m_2$ between morphisms directly.\foot{Since $m_2$ has degree
zero with respect to the ${\bf Z}$-grading on morphisms, we
can ignore the universal cover of the space of lagrangian
planes and just consider the slope of the lines.}
For the moment we set all translations (which parametrize
bundles of the same rank and degree) to zero, which means
in part, on the Fukaya side (as we shall see)
that all lines pass through the origin (we assume an origin has
been chosen for the elliptic curve and its mirror).
We also set all holonomies to be zero (trivial local systems).
We will return to
non-zero translations and holonomies later in this section.
We define the K\"ahler parameter to be
$$\rho = b + iA,$$ where $A$ is the area of the torus,
and $b$ describes a two-form of the same name:  $b = b\, dx\wedge dy \in
H^2(\widetilde{E};{\bf R})/H^2(\widetilde{E};{\bf Z}).$
In terms of the complexified K\"ahler form $\omega$ of the
previous section, $\rho = \int_{\widetilde{E}}\omega$ (and $A = k$).
We set $b=0$ to begin with, so $\rho = iA.$
The mirror map says that $\widetilde{E}$ is dual to $E$ when
$\rho = \tau.$
  
Let
$$\eqalign{{\cal L}_1 &= (1,0),\cr {\cal L}_2 &= (1,1),\cr
{\cal L}_3 &= (1,2).}$$
On the derived category side, then, we are considering
line bundles of degrees $0,1,2.$
Note that ${\cal L}_1 = {\cal O},$ the sheaf of holomorphic
functions.
We can define $L \equiv {\cal L}_2;$ then ${\cal L}_3 = L^2,$
and
$$\eqalign{\hom({\cal L}_1,{\cal L}_2) &= H^0(L),\cr
\hom({\cal L}_2,{\cal L}_3) &= H^0(L),\cr
\hom({\cal L}_1,{\cal L}_3) &= H^0(L^2).}$$
The product of global sections gives us a map
\eqn\mtwo{m_2:  H^0(L)\otimes H^0(L) \longrightarrow H^0(L^2).}
The theta function $\theta(\tau,z)$ is the unique global
section of $L,$ and the decomposition of the product of theta
functions into two sections of $L^2$ is known as the
``addition formula.''
  
Let's see what we get from the Lagrangian lines ${\cal L}_i$
(we use the same notation as for bundles)
We have three lines through the origin of slope $0, 1,$ and $2$
respectively.  Note that ${\cal L}_1\cap {\cal L}_2 = \{{\rm e}_1\},$
where the origin ${\rm e}_1 \equiv (0,0).$  Also, ${\cal L}_2\cap {\cal L}_3 = \{{\rm e}_1\},$
while
${\cal L}_2\cap {\cal L}_3 = \{{\rm e}_1,{\rm e}_2\},$
where ${\rm e}_2 \equiv (1/2,0).$
On the left hand side of \mtwo,
${\rm e}_1$ represents the theta function $\theta(\tau,z)$.\foot{As
we point out again in the next section, the notation is ambiguous:
the section being represented by ${\rm e}_1$ depends on the two
lines of which we take it as an intersection point.}
On the right hand side, the ${\rm e}_i$ represent
a distinguished basis of the two-dimensional space of
sections of $L^2$ (see section {\sl 2.3} for details).
We define this basis correspondence by
$$\eqalign{{\rm e}_1 &\leftrightarrow \theta[0,0](2\tau,2z),\cr
{\rm e}_2 &\leftrightarrow \theta[1/2,0](2\tau,2z).}$$
We will use the mirror map $$\tau = \rho$$ to check the correspondence
of the $m_2$ product.
  
The map $m_2$
is:
$$m_2({\rm e}_1\otimes{\rm e}_1) = C({\rm e}_1,{\rm e}_1,{\rm e}_1)\cdot {\rm e}_1
+ C({\rm e}_1,{\rm e}_1,{\rm e}_2)\cdot {\rm e}_2,$$
where we use Fukaya's procedure \mcoef\ to compute the matrix elements, $C.$
The holomorphic maps from the discs are specified by the triangles
bounded by ${\cal L}_1,{\cal L}_2,$ and ${\cal L}_3$ with vertices given
by the arguments of $C.$  Thus, $C({\rm e}_1,{\rm e}_1,{\rm e}_1)$
is given by summing all triangles with the origin as vertex.  Looking
at the universal cover of the torus, we need (up to translation) triangles
with lattice points as vertices and sides of slope $0, 1,$ and $2.$
Consider the base of the triangle, which must be integral length,
say $n,$ running from $(0,0)$ to $(n,0).$  The third vertex is
$(2n,2n),$ which is a lattice point.  Thus all the triangles are
indexed by the integers, and the $n^{th}$ triangle has area
$n^2$ times the area $A.$  To determine $C({\rm e}_1,{\rm e}_1,{\rm e}_2)$
we use the triangles with vertices $(0,0),(0,n+1/2),(2n+1,2n+1)$
and area $A(n + 1/2)^2.$
  
Therefore,
$$\eqalign{C({\rm e}_1,{\rm e}_1,{\rm e}_1) &= \sum_{n=-\infty}^{\infty}
\exp[-2\pi An^2]\cr
C({\rm e}_1,{\rm e}_1,{\rm e}_2) &= \sum_{n=-\infty}^{\infty}
\exp[-2\pi A(n+1/2)^2].}$$
Note that the triangles for $n\geq 0$ and $n<0$ are related by
the ${\bf Z}_2$ automorphism of the system (torus plus configuration
of Lagrangian lines); nevertheless, we treat them as distinct.
We recognize (see section {\sl 2.3}) the coefficients $C$ as theta functions
evaluated at $z = 0$ (with no shifts).  Specifically,
$$\eqalign{C({\rm e}_1,{\rm e}_1,{\rm e}_1) &=
\theta[0,0](i2A,0) = \theta[0,0](2\rho,0),\cr
C({\rm e}_1,{\rm e}_1,{\rm e}_2) &= \theta[1/2,0](i2A,0) =
\theta[1/2,0](2\rho,0).}$$
The product $m_2$ therefore precisely contains the information
of the addition formula
$$\theta(\rho,z)\theta(\rho,z) = \theta[0,0](2\rho,0)\theta[0,0](2\rho,2z) +
\theta[1/2,0](2\rho,0)\theta[1/2,0](2\rho,2z).$$
  
We now consider the case where the third line doesn't pass through the
point of intersection of the other two.  Call $\alpha$ the
(positive) $x$-intecept
of the closest line of slope two, so
$0<\alpha < 1/2.$  Then the points
of ${\cal L}_1 \cap{\cal L}_3$ are $(\alpha, 0),(1/2 + \alpha,0).$
If we label the $\alpha$-shifted basis ${\rm e}_{1,\alpha},
{\rm e}_{2,\alpha},$ then we have
$$\eqalign{C({\rm e}_1,{\rm e}_1,{\rm e}_{1,\alpha}) &= \sum_{n=-\infty}^{\infty}
\exp(-2\pi A[n+\alpha)^2] = \theta[\alpha,0](2\rho,0)\cr
C({\rm e}_1,{\rm e}_1,{\rm e}_{1,\alpha}) &= \sum_{n=-\infty}^{\infty}
\exp[-2\pi A(n+1/2 + \alpha)^2] =
\theta[1/2 + \alpha,0](2\rho,0).}$$
One also needs to know the rules for identifying the sections
corresonding to the shifted basis, and the obvious guess is
$$\eqalign{{\rm e}_{1,\alpha}&\leftrightarrow \theta[\alpha,0](2\tau,2z)\cr
{\rm e}_{2,\alpha}&\leftrightarrow\theta[1/2 + \alpha,0](2\tau,2z).}$$
In fact, this is almost correct.  Some phases need to be added, after
which the composition becomes once again equivalent to the addition
formula, now with $\alpha$-dependent shifts, as indicated.
The exact formulas are given in the next section.
The appearance of these phases is not unexpected, since
there is a certain arbitrariness in the choice of $\tau$-dependence
of the phase in the definition of the theta function.
  
Let us return to zero translation, but add holonomy
to ${\cal L}_3.$  That is, we take the flat line bundle with
connection $2\pi i \beta dt_3$ where $t_3 \sim t_3 + 1$ is a coordinate
along ${\cal L}_3.$  The matrix elements $C$ then have to be
weighted by $\exp[-2\pi ({\rm Area}) + 2\pi i\oint \beta dt].$
This gives
$$\eqalign{C({\rm e}_1,{\rm e}_1,{\rm e}_1) &= \sum_{n=-\infty}^{\infty}
\exp[-2\pi A(n+\alpha)^2 + 2\pi in\beta]
= \theta[0,\beta](2\rho,0),\cr
C({\rm e}_1,{\rm e}_1,{\rm e}_2) &= \sum_{n=-\infty}^{\infty}
\exp[-2\pi A(n+\alpha)^2 + 2\pi i(n+1/2)\beta]
= \theta[1/2,\beta](2\rho,0).}$$
The basis correspondence
the involves $\theta[0,\beta](2\tau,2z)$ and $\theta[1/2,\beta](2\tau,2z).$
Now it is clear that the translated case with holonomy gives
$$\eqalign{
C({\rm e}_1,{\rm e}_1,{\rm e}_{1,\alpha}) &=
\theta[\alpha,\beta](2\rho,0),\cr
C({\rm e}_1,{\rm e}_1,{\rm e}_{2,\alpha}) &=
\theta[1/2 + \alpha,\beta](2\rho,0),}$$
and that the ${\rm e}_{i,\alpha}$ correspond to $\theta[\alpha,\beta]$
and $\theta[1/2 + \alpha,\beta]$
(up to phases).  But these theta functions are
precisely those for the line bundles of degree two described by
the point $\alpha\tau + \beta$ on the Jacobian torus.
We therefore have learned how the shift functor representing
the moduli of different bundles of same rank and degree (i.e. the
Jacobian) acts on the Fukaya side by shifting lines, adding
holonomy, and introducing phases.
  
We see that the two real translation parameters describing
the Jacobian torus, which are on the same footing from the
point of view of line bundles, have different interpretations
on the Lagrangian side -- in terms of holonomies and displacements
of the Lagrangian submanifold.  This is perhaps not too surprising
if we consider the string theoretic origin of the moduli of the
D-brane.  Both moduli come from zero modes associated to the
ten-dimensional gauge field.  On the D-brane the transverse components
become normal vectors and the zero modes give the motions of the brane,
while the parallel components describe the flat bundle moduli (holonomy).
  
It is now a simple matter to check that these formulas are true for
general $\rho$ (recall $\rho = \int_{\widetilde{E}}\omega$),
i.e. when $b \neq 0$, since we weight a holomorphic map from the
disc by $\exp[2\pi i\int \phi^*\omega + \oint \beta dt].$
Complexifying the K\"ahler class in this way is familiar from mirror
symmetry, and the holonomy is expected for D-branes or open strings.
We now have a clear description of all the parameters involved
in the mirror map $\tau \leftrightarrow \rho$ as well as the translational bundle (or Lagrangian)
moduli.
  
The composition is somewhat
more involved for stable vector bundles of rank $r$ and degree $d$
(which correspond to lines $(r,d)$ with
${\rm gcd}(r,d) = 1$), but as we will see in the
next section, these cases can be deduced from our knowledge of line
bundles.  The non-stable bundles
(which correspond to non-unitary local systems) are more subtle.
Identifying a proper basis for the global sections is the
main difficulty.  Fortunately, we will be able to describe all
sections in terms of simple theta
functions.  The Fukaya compositions are rather easy to compute, as above,
and the equivalence of these products
to the decomposition of sections of bundles will be shown
to reduce to the classical addition formulas, as above.
  
\newsec{Categorical Equivalence for the Elliptic Curve}
  
We must construct the categorical equivalence in the general case,
where the objects in the derived category are more general than
line  bundles.  We can immediately reduce to indecomposable
bundles, by linearity of $m_2.$
The  idea, then,  is  to  use the representation of
indecomposable vector
bundles on elliptic curve as push-forwards under the isogenies
of line bundles tensored with bundles of type $F(V,\exp(N))$ where
$N$ is a constant nilpotent matrix (see section {\sl 2.3}).
Let us denote by ${\cal L}(E)$ the category of bundles on $E$ of the form
$L(\phi)\otimes F(V,\exp(N))$.
We claim that it is sufficient to construct a natural
equivalence of ${\cal L}(E)$ with the appropriate
subcategory in the Fukaya category. More precisely, we
need that these equivalences for elliptic curves
$E_{q^r}$ and $E_q$ commute with the corresponding
pull-back functors $\pi_r^*$.
 
{\bf Main Theorem}
 
The categories
${\cal D}^b(E_q)$ and ${\cal F}^0(\widetilde{E}^q)$ are equivalent:
$${\cal D}^b(E_q)\cong{\cal F}^0(\widetilde{E}^q)$$
(we use $q = \exp(2\pi i \tau)$ instead of $\tau$ here).
 
We prove this by constructing a bijective functor
$\Phi: {\cal D}^b(E_q)\rightarrow{\cal F}^0(\widetilde{E}^q)$
in five steps.  We first define it on a restricted class of objects:
bundles of the form $L(\varphi)\otimes F(V,\exp N)$ (step one).  Then we show
how this respects the composition law (step two).  Step three is showing
that it commutes with isogenies.  Then, in step four we treat general
vector bundles by constructing them from isogenies.  Finally, in step
five, we treat
the only remaining case:  torsion sheaves, or
thickened skyscrapers (recall that objects in the derived
category are, up to
shifts, direct sums of vector bundles and these sheaves, which have
support only at points).
 
More specifically,
we construct a fully faithful functor from the abelian
category of coherent sheaves to the Fukaya category.  The
choice of logarithms of slopes for the corresponding objects
is $\pi i \alpha,$ where $\alpha \in (-1/2,1/2]$
(the right boundary is achieved by torsion sheaves).
On this subcategory $\hom^0$ between lines with slopes
$\exp(\pi i \alpha_1)$ and $\exp(\pi i \alpha_2)$ is
non-zero only if $\alpha_1 < \alpha_2.$  This functor can
be extended to the equivalence of the entire derived
category of coherent sheaves with the Fukaya category, using
Serre duality and the rule that the shift by one in the derived
category corresponds to $\alpha \mapsto \alpha + 1.$
 
\subsec{The Functor $\Phi$}
 
Let ${\cal L}(E_q)$ denote the full subcategory of the category of vector
bundles on $E_q$ consisting of bundles of the form
$L(\varphi)\otimes F(V,\exp N),$ where
$\varphi = t_x^*\varphi_0\cdot\varphi_0^{n-1}$ for some $x\in E_q,
n\in \bfz$ (we use the conventions of section {\sl 2.3}.)
 
We construct a functor
$$\Phi_q: {\cal L}(E_q) \rightarrow {\cal F}^0(\widetilde{E}^q),$$
where ${\cal F}^0$ denotes the Fukaya category and $q$ denotes the
exponential of the complexified K\"ahler parameter on the right
hand side ($q = \exp(2\pi i \rho),$ and $\rho = b + ik$ is set equal
to $\tau,$ as per the mirror map).  On the Fukaya side,
an object is a pair of line $\Lambda$ in $\bfr^2/\bfz^2$
described by a parametrization $(x(t),y(t))$
in $\bfr^2$ and a flat connection $A$ on the line.  We describe the
flat connection as a constant, $V$-valued one-form on $\bfr^2$ --
the restriction to $\Lambda$ is implied.
Specifically, the map $\Phi$ on objects is
$$\eqalign{
\Phi:  \left(L(t_{\alpha \tau+\beta}^*\varphi_0\cdot\varphi_0^{n-1})
\otimes F(V,\exp N)\right) \mapsto
(\Lambda, A),\qquad \Lambda &= (\alpha + t, (n-1)\alpha + nt),\cr
A &= (-2\pi i \beta \cdot {\bf 1}_V + N) dx}$$
(we subsequently drop the notation ${\bf 1}_V$ for the identity map on $V$).
In other words, the line has slope $n$ (degree of the line bundle)
and $x$-intercept $\alpha/n.$  The monodromy matrix between points
$(x_1, y_1)$ and $(x_2,y_2)$ is
$\exp[-2\pi i \beta (x_2-x_1)]\exp[N(x_2-x_1)].$  This is
precislely what we saw last section: the shifts are represented
as translations of the line and as monodromies.
 
It remains to describe the map on morphisms:
$$\eqalign{\Phi:  \hom\big( L(\varphi_1)\otimes &F(V_1,\exp N_1),
L(\varphi_2)\otimes F(V_2,\exp N_2)\big)\rightarrow \cr
&
\hom^0\big(
\Phi(L(\varphi_1)\otimes F(V_1,\exp N_1)),\Phi(L(\varphi_2)
\otimes F(V_2,\exp N_2))\big),}$$
where $\varphi_i = t_{\alpha_i \tau + \beta_i}^*\varphi_0
\cdot\varphi^{n_i-1},$
$i = 1,2,$
and $\hom^0$ is the image in ${\cal F}^0,$ i.e. the degree zero
part of $\hom.$  Now recall Prop. 2.
If $n_1>n_2,$ both of the above spaces are zero.
If $n_1 = n_2,$ then either the above spaces are zero, or
$L(\varphi_1)\cong L(\varphi_2)$
and the problem reduces to homomorphisms of vector spaces.
If $n_1<n_2,$ then
$$\eqalign{LHS &= H^0\left(
L(\varphi_2\varphi_1^{-1})\otimes
F(V_1^*\otimes V_2, \exp(N_2-N_1^*))\right)\cr
&= H^0(L(\varphi_2\varphi_1^{-1}))\otimes (V_1^*\otimes V_2)\qquad {\rm by }\;
{\cal V},}$$
while
$$RHS = \bigoplus_{{\rm e}_k\in \Lambda_1\cap \Lambda_2} V_1^*\otimes V_2
\; \cdot {\rm e}_k.$$
The points of $\Lambda_1 \cap \Lambda_2$
are easily found from $\Phi$ to be
$${\rm e}_k = \left( {k + \alpha_2 - \alpha_1 \over n_2 - n_1},
{n_1k + n_1\alpha_2 - n_2\alpha_1 \over n_2 - n_1}\right),
\qquad k \in \bfz/(n_2-n_1)\bfz.$$
We note that
$$\varphi_2 \varphi_1^{-1}=
t_{\alpha_2 \tau + \beta_2}^*\varphi_0\cdot \varphi_0^{n_2-1}
\cdot t_{\alpha_1\tau+\beta_1}^*\varphi_0^{-1}\varphi_0^{-n_1+1}
=t_{\alpha_{12} \tau + \beta_{12}}^*\varphi_0^{n_2-n_1},$$
where
$$\alpha_{12} = {\alpha_2 - \alpha_1\over n_2-n_1},\qquad
\beta_{12} = {\beta_2 - \beta_1 \over n_2 - n_1}.$$
Now we have the standard basis of theta functions on $H^0(L(\varphi_2\varphi_1^{-1})):$
$$\eqalign{t_{\alpha_{12}\tau + \beta_{12}}^*
\theta\left[{k\over n_2-n_1},0\right]&\big((n_2-n_1)\tau,(n_2-n_1)z\big)
= \cr
&\theta\left[{k\over n_2-n_1},0\right]
\big( (n_2-n_1)\tau,(n_2-n_1)(z+ \alpha_{12}\tau + \beta_{12})\big),}$$
$k\in \bfz/(n_2-n_1)\bfz.$
Let us call this function $f_k.$  The standard basis on the Fukaya
side is also indexed by $k.$  We recall that $\alpha_{12}$ and $\beta_{12}$
effect shifts and monodromies on the right hand side, and this determines
the identification of bases up to a constant.  Let $T \in V_1^*\otimes V_2.$
We define
$$\Phi\left( {\cal V}(f_k\otimes T)\right)
= \exp(-\pi i \tau \alpha_{12}^2(n_2-n_1))\exp[\alpha_{12}(
N_2-N_1^*- 2\pi i (n_2-n_1)\beta_{12})]\cdot T\; {\rm e}_k.$$
 
\subsec{$\Phi\circ m_2 = m_2 \circ \Phi$}
 
We must check that the definition of $\Phi$ respects the composition
maps in the two categories.
We have $L(\varphi_i)\otimes F(V_i,\exp N_i),$ $i = 1...3.$
 
a)  We have to compare compositions of sections
$${\cal V}\left( t_{\alpha_{12}\tau + \beta_{12}}^*
\theta\left[{a\over n_2-n_1},0\right]
\big((n_2-n_1)\tau\, ,\, (n_2-n_1)z\big)\otimes A
\right)$$
and
$${\cal V}\left( t_{\alpha_{23}\tau + \beta_{23}}^*
\theta\left[{a\over n_3-n_2},0\right]
\big((n_3-n_2)\tau\, ,\, (n_3-n_2)z\big)\otimes B
\right),$$
where $A \in V_1^*\otimes V_2$ and
$B\in V_2^* \otimes V_3.$
Using Prop. 3, this composition is equal to ${\cal V}$ applied to the
following expression:
$$\eqalign{
{\rm Tr}_{V_2}\theta\Big[&{a\over n_2-n_1},0\Big]\left((n_2-n_1)\tau\, ,\, (n_2-n_1)
(z+\alpha_{12}\tau + \beta_{12} - {N\over n_2-n_1})\right)\times \cr
&\phantom{X}
\times
\theta\Big[{b\over n_3-n_2},0\Big]\big((n_3-n_2)\tau\, ,\, (n_3-n_2)
(z+\alpha_{23}\tau + \beta_{23} - {N\over n_3-n_2})\big)\cdot
 A\otimes B,}
$$
where
$$\eqalign{N &= {1\over 2\pi i}{ -(n_2-n_1)(N_3 - N_2) + (n_3-n_2)(N_2-N_1^*)\over n_3-
n_1}
\cr &= {1\over 2\pi i}\left(N_2 - {(n_2-n_1)N_3 + (n_3-n_2)N_1^*\over n_3-n_1}\right)}$$
(we have replaced $N_2^*$ by $N_2$ under the trace sign
using $\langle v,N_2^* \xi\rangle = \langle N_2 v,\xi\rangle$ for
$v\in V_2,\xi \in V_2^*$).
Now the addition formula (II.6.4 of \mum)
implies the following $\theta$-identity:
$$\eqalign{
\theta\left[{a\over n_2-n_1},0\right]&
\big((n_2-n_1)\tau\, ,\, (n_2-n_1)z_1)\big)\cdot
\theta\left[{b\over n_3-n_2},0\right]
\left((n_3-n_2)\tau\, ,\, (n_3-n_2)z_2\right) =
\cr
&=
\sum_{m\in \bfz}\exp
\left[\pi i \tau{k_m^2\over (n_2-n_1)(n_3-n_2)(n_3-n_1)}
+ 2\pi i {k_m\over n_3-n_1}(z_2-z_1)\right]
\times \cr
&
\phantom{\sum_{m\in \bfz},}
\theta\left[{a+b+(n_3-n_2)m\over n_3-n_1},0\right]
\big( (n_3-n_1)\tau\, ,\, (n_2-n_1)z_1 + (n_3-n_2)z_2\big),
}$$
where $$k_m = (n_2-n_1)b - (n_3-n_2)a + (n_2-n_1)(n_3-n_2)m.$$
Plugging in $z_1 = z+\alpha_{12}\tau + \beta_{12} - {N\over n_2-n_1}$
and $z_2 = z+\alpha_{12}\tau + \beta_{12} - {N\over n_2-n_1},$
we obtain
$$\eqalign{
{\rm Tr}_{V_2}
\sum_{m\in \bfz}\exp
\Big[&{(\pi i \tau) k_m^2\over (n_2-n_1)(n_3-n_2)(n_3-n_1)}\cr
&+ {(2\pi i) k_m\over n_3-n_1}\big( (\alpha_{23}-\alpha_{12})\tau
+ \beta_{23}-\beta_{12} + {n_3-n_1\over (n_2-n_1)(n_3-n_2)}N\big)\Big]
\cdot A\otimes B \times \cr
&\times\theta\left[ {a+b+(n_3-n_2)m\over n_3-n_1},0\right]
\big( (n_3-n_1)(z + \alpha_{13}\tau + \beta_{13})\big).}$$
 
b)  In the Fukaya category, we need to compose $A \cdot {\rm e}_a$
and $B \cdot {\rm e}_b,$ where
$A \in V_1^*\otimes V_2$ and $B\in V_2^*\otimes V_3$
and $a\in \bfz/(n_2-n_1)\bfz$ and
$b\in \bfz/(n_3-n_2)\bfz$
label points in $\Lambda_1\cap \Lambda_2$
and $\Lambda_2 \cap \Lambda_3,$ respectively.\foot{There
is some notational ambiguity:  ``What is ${\rm e}_3$?''
The context should make it clear which intersections are being
indexed, and therefore what the range of the index should be.}
 
The points of intersection are easily computed.
For two lines $\Lambda_i$ and $\Lambda_j,$ the
$n_j-n_i$ points of intersection in $\bfr^2/\bfz^2$
are
$${\rm e}_k = \left(\alpha_{ij} + {k\over n_j-n_i},
{n_i\alpha_j-n_j\alpha_i+ n_ik \over n_j-n_i}\right).$$
We must sum over triangles in the plane formed by the
$\bfz$-translates of $\Lambda_3$ which pass through points
lattice-equivalent to ${\rm e}_b.$
Note  that  such  a  triangle  is uniquely determined by the
first coordinate of the intersection point of $\Lambda_3$
with $\Lambda_2$ which should be of the form
$\alpha_{23}+{b\over n_3-n_2}+m$, $m\in \bfz.$
The point of intersection
${\rm e}_c$ in $\Lambda_1\cap \Lambda_3$ depends on $m,$ and
is calculated to be given by the above formula with
$c = a + b + m(n_3-n_2) \;{\rm mod }\;(n_3-n_1).$
 
If we call $l_1$ the horizontal distance (difference of $x$ coordinates)
between ${\rm e}_a$ and ${\rm e}_c,$
$l_2$ the horizontal distance between ${\rm e}_a$ and ${\rm e}_b,$
and $l_3$ the horizontal distance between ${\rm e}_b$ and ${\rm e}_c$
(may be negative), then by adding up $x$ and $y$ differences we have
$$\eqalign{l_1 &= l_2 + l_3\cr n_1 l_1 &= n_2 l_2 + n_3 l_3}.$$
Therefore, $l_1 = {n_3-n_2\over n_3-n_1}l_2,$
and $l_3 = {n_1-n_2\over n_3-n_1}l_2.$  Now we know $l_2$ explicitly
from above.  The area $\Delta$ of the triangle (rather, the number of
fundamental domains) formed by the vertices is
half that of the parallelogram determined by vectors
$(l_1,n_1l_1)$ and $(l_2,n_2l_2).$
Hence,
$$\Delta = {1\over 2}{\rm Det}\pmatrix{l_1&l_2\cr n_1l_1&n_2 l_2}
= l_1l_2(n_2-n_1)/2 = {(n_3-n_1)(n_2-n_1)\over 2(n_3-n_2)} (l_2)^2.$$
Plugging in $l_2 = \alpha_{23}-\alpha_{12} + b/(n_3-n_2) - a/(n_2 - n_1)
+ m,$ gives
$$\Delta = {1\over 2}{\left[ k_m + (\alpha_{23}-
\alpha_{12})(n_3-n_2)(n_2-n_1)\right]^2\over
(n_3-n_2)(n_3-n_1)(n_2-n_1)},$$
where $k_m$ was defined above.
Now since $\int \phi^*\omega = \rho \Delta = \tau \Delta,$
then recalling the assigned monodromies, we have,
by \mcoef\ for the Fukaya product
(recalling the sign conventions of the $l_i$),
$$\eqalign{A \cdot {\rm e}_a \circ B \cdot {\rm e}_b &= \cr
\sum_{m\in\bfz}\exp &\lbrace\pi i \tau
{\left[ k_m + (\alpha_{23}-
\alpha_{12})(n_3-n_2)(n_2-n_1)\right]^2\over
(n_3-n_2)(n_3-n_1)(n_2-n_1)}\rbrace\times
\cr
&{\rm Tr}_{V_2}
\exp\left[
-l_1(N_1^* - 2\pi i \beta_1)
+ l_2 (N_2 - 2\pi i \beta_2)
+ l_3 (N_3 - 2\pi i \beta_3)\right]\cdot\cr
&\phantom{{\rm Tr}_{V_2}}\cdot A\otimes B \cdot
{\rm e}_{a + b + m(n_3-n_2)}.}$$
Upon expanding the square,
substituting the correct values of the $l_i,$
and recalling the definition of $N,$ the right hand side becomes
$$\eqalign{
{\rm Tr}_{V_2}
\sum_{m\in \bfz}
\exp &\Big[
\pi i \tau {k_m^2\over (n_3-n_2)(n_3-n_1)(n_2-n_1)}
+ 2\pi i {k_m\over n_3-n_1}\big( (\alpha_{23}-\alpha_{12})\tau
+ \beta_{23}-\beta_{12}\big)\cr
&\phantom{\big[}+ \pi i \tau
(\alpha_{23}-\alpha_{12})^2(n_3-n_2){n_2-n_1\over n_3-n_1}\cr
&+ 2\pi i (\alpha_{23}-\alpha_{12})(\beta_{23}-\beta_{12})
{(n_3-n_2)(n_2-n_1)\over n_3-n_1} \cr
&\phantom{\Big[}
+ {k_m + (\alpha_{23}-\alpha_{12})(n_3-n_2)(n_2-n_1)\over
(n_3-n_2)(n_2-n_1)}2\pi i N \big]\cdot A\otimes B\cdot
{\rm e}_{a + b + m(n_3-n_2)}.}$$
Taking into account the exponential factors in the definition of
$\Phi$ yields equivalence with the product in the derived category,
which is the desired result.
 
\subsec{Isogeny}
 
On our subcategory ${\cal L},$ we have the functor of pull-back under
isogeny:
$$\pi_r^*: {\cal L}(E_q) \rightarrow {\cal L}(E_{q^r}).$$
We describe the analogue on the Fukaya side.
Consider $\pi_r: \bfr^2/\bfz^2\rightarrow \bfr^2/\bfz^2,$
which sends $(x,y)$ to $(rx,y).$  This is an $r$-fold covering
and respects the K\"ahler form if the complexified K\"ahler
parameter on the left hand side is $r\rho$ when that of the
right hand side is $\rho.$
We then have the corresponding functors $\pi_r^*$ and
$\pi_{r*}$ between Fukaya categories.
 
{\bf Proposition 4}
 
One has the following commutative diagram of
functors:
$$\matrix{{\cal L}(E_{q})&\matrix{ {}_\Phi \cr \longrightarrow \cr {}}&
{\cal F}(\widetilde{E}^\rho)\cr
\matrix{\pi_r^* &\downarrow & {}}&{}&\matrix{{}&\downarrow&\pi_r^*}\cr
{\cal L}(E_{q^r})&\matrix{ {}_\Phi \cr \longrightarrow \cr {}}&
{\cal F}(\widetilde{E}^{r\rho}).}$$
 
Proof:  This is a straightforward check using the formula
$$\theta[a/n,0](n\tau,nz) = \sum_{k\in \bfz/r\bfz}
\theta[(a+nk)/(nr),0](nr^2\tau,nrz).$$
 
\subsec{The General Case}
 
Now we can extend the functor $\Phi$ from the category ${\cal L}(E_q)$
to the category of all vector bundles.  For this we use the fact that
any indecomposable bundle on $E_q$ is isomorphic to
$$\pi_{r*}\left( L_{q^r}(\varphi)\otimes F_{q^r}(V,\exp N)\right)$$
(see Prop. 1).
Hence we can set
$$\Phi_q\left( \pi_{r*}[L_{q^r}(\varphi)\otimes F_{q^r}(V,\exp N)]\right)
= \pi_{r*}\left( \Phi_{q^r}[L_{q^r}(\varphi)\otimes
F_{q^r}(V,\exp N)]\right).$$
 
It remains to define $\Phi_{q}$ on morphisms between vector bundles.
For this we use the isomorphism \homeq:
$$\eqalign{\hom&\left(
\pi_{{r_1}*}[L_{q^{r_1}}(\varphi)\otimes F_{q^{r_1}}(V_1,\exp N_1)],
\pi_{{r_2}*}[L_{q^{r_2}}(\varphi)\otimes F_{q^{r_2}}(V_2,\exp N_2)]\right)
\cong \cr
&\hom\left(
\pi_{{2}}^*[L_{q^{r_1}}(\varphi)\otimes F_{q^{r_1}}(V_1,\exp N_1)],
\pi_{{1}}^*[L_{q^{r_2}}(\varphi)\otimes F_{q^{r_2}}(V_2,\exp N_2)]\right),}$$
where $\pi_1$ and $\pi_2$ are defined from the cartesian square
$$\matrix{E_{12}&\matrix{{}_{\pi_1}\cr\longrightarrow\cr{}}&E_{q^{r_2}}\cr
\matrix{\pi_2&\downarrow&{}}&{}&\downarrow\cr
E_{q^{r_1}}&\longrightarrow &E_q.}$$
Here $E_{12}$ is an elliptic curve if and only if ${\rm gcd}(r_1,r_2) = 1.$
In general, $E_{12} \cong E_{q^{r_1'r_2'}}\times \bfz/d\bfz,$
where $d = {\rm gcd}(r_1,r_2)$ and $r_i' = r_i/d,$ $i = 1,2.$
Thus, the above Hom space is decomposed into a direct sum of
$d$ Hom spaces between objects of the category ${\cal L}(E_{q^{r_1'r_2'}}).$
There is a similar decomposition of the corresponding Hom spaces
on the Fukaya side, so we just take the direct sum of isomorphisms
between Hom's given by $\Phi$ on ${\cal L}(E_{q^{r_1'r_2'}}).$
To check that this extended map $\Phi$ is still compatible with
compositions we note that given a triple of
vector bundles of the form
$\pi_{r_i*}\left(L_{q^{r_i}}(\varphi_i)\otimes
F_{q^{r_i}}(V_i,\exp N_i)\right),$
$i = 1...3,$ we can embed the pairwise Hom's into the corresponding
Hom spaces between the pull-backs of
$L_{q^{r_i}}(\varphi_i)\otimes
F_{q^{r_i}}(V_i,\exp N_i)$ to the triple fibered product of $E_{q^{r_1}},
E_{q^{r_2}},$ and $E_{q^{r_3}}$ over $E_q,$ which is still a disjoint
union of elliptic curves.  Now the required compatibility follows from
the fact that the functor $\Phi$ on the ${\cal L}(E_q)$ commutes with
the pull-back functor $\pi_r^*.$
 
\subsec{Extension to Torsion Sheaves}
 
For every $z_0 \in \bfc$ and a nilpotent operator
$N \in {\rm End}(V)$ we have the corresponding coherent sheaf of $\bfc$
supported at $z_0.$  Namely, ${\cal O}_{rz_0}\otimes V/\langle
z - z_0 - {N\over
2\pi i}\rangle,$
where $r = {\rm dim}V$ is the smallest positive integer such that $N^r = 0.$
For example, if $N=0$ then the notation means that we set to zero
anything in the
ideal generated by $z - z_0,$ i.e. we get ${\cal O}_{z_0}$ as in \sky.
We denote by $S(z_0,V,N)$ the
direct image of this sheaf on $E.$  Then for every object
$L(\varphi)\otimes F(V',\exp N')$ of ${\cal L}(E)$
we have canonically
$\hom\big( L(\varphi)\otimes F(V',\exp N'), S(z_0,V,N)\big)
\cong \hom(V',V).$
The composition map
$$\eqalign{\hom &\left( L(\varphi_1)\otimes F(V_1,\exp N_1),
L(\varphi_2)\otimes F(V_2,\exp N_2)\right)\otimes \cr
&\phantom{\hom()\hom()}\hom\left( L(\varphi_2)\otimes F(V_2,\exp N_2),S(x,V,N)\right)
\longrightarrow\cr
&\phantom{\big(L(\varphi_1)\otimes F(V_1,\exp N_1),\big)}
\hom\left( L(\varphi_1)\otimes F(V_1,\exp N_1),
S(x,V,N)\right)}$$
can be written (using the isomorphism ${\cal V}$) as follows:
$${\cal V}(f\otimes A)\circ B = {\rm Tr}_{V_2}
f\left(x + {N\over 2\pi i} - {1\over 2\pi i (n_2-n_1)}(N_2-N_1)^*\right)
(A\otimes B),$$
where $f\in H^0(L(\varphi_2\varphi_1^{-1})),$ $A\in \hom(V_1,V_2),$
and $B\in \hom(V_2,V).$
 
We extend $\Phi$ to torsion sheaves by sending $S(\alpha\tau + \beta,V,N)$
to the line $(-\alpha,t)$ with connection given by 
$(2\pi i \beta \cdot {\bf 1}_V + N)dy.$
We also define an isomorphism
$$\eqalign{\Phi: \hom &\left(
L(\varphi_1)\otimes F(V_1,\exp N_1),S(x,V_2,N_2)\right)\longrightarrow\cr
&\hom\left( \Phi(L(\varphi_1)\otimes F(V_1,\exp N_1)),
\Phi(S(x,V_2,N_2))\right) \cong \hom(V_1,V_2)}$$
by the formula (where $\varphi_1 = t_{\alpha_1 \tau +
\beta_1}^*\varphi_0\cdot\varphi_0^{n-1}$
and $x = \alpha_2\tau + \beta_2$)
$$\eqalign{A \mapsto
\exp[-\pi i \tau n \alpha_2^2 - 2\pi i\tau \alpha_1\alpha_2
- &\alpha_2(nN_2-N_1^*+2\pi i (\beta_1+n\beta_2)) - \cr 
&\alpha_1 (N_2+2\pi i\beta_2)]\cdot A,}$$
where $A\in \hom(V_1,V_2).$
 
It remains to check that $\Phi$ respects compositions.  We only have
to check this for the composition of $\hom$'s between
$L(\varphi_1)\otimes F(V_1,\exp N_1),$
$L(\varphi_2)\otimes F(V_2,\exp N_2),$
and $S(x,V,N).$  The proof of this is similar to the cases
considered previously, but even simpler.  One needs only
the definition of the theta function; no theta identities are
used.  We leave the case of $\hom$'s from higher rank vector
bundles as an exercise for the reader (use $\pi_{r*}$).
 
{\bf Remark}
 
Any automorphism $\varphi$ of $\mir$ preserving
the K\"ahler structure and $\Omega$ induces an
autoequivalence $\varphi_*$ of the Fukaya category
$\fuk$ defined up to a shift by an integer.  In the case
of the flat torus $T = \bfr^2/\bfz^2,$ this leads to an
action of the central extension of $SL(2;\bfz)$ by
$\bfz$ in the Fukaya category ${\cal F}(T).$  The corresponding projective
action of $SL(2;\bfz)$ on ${\cal D}^b(E)$ is given as follows.
The matrix $\pmatrix{1&0\cr 1 &1}$ acts by tensoring with $L,$
while the matrix $\pmatrix{0 & 1\cr -1 &0}$ acts as the
Fourier-Mukai transform
$${\cal S}\mapsto p_{2*}(p_1^*{\cal S}\otimes {\cal P}),$$
where ${\cal P}$ is the normalized Poincar\'e line bundle on
$E\times E$ given by
$${\cal P} = m^*{\cal O}(-x_0)\otimes p_1^*{\cal O}(x_0)
\otimes p_2^*{\cal O}(x_0),$$
where $m: E\times E\rightarrow E$ is the group law,
$x_0$ is the point corresponding to $1/2 + \tau/2,$ and
$p_1, p_2$ are the projections.

\newsec{Conclusions}
  
We have shown that for the elliptic curve,
as Kontsevich conjectured, mirror symmetry has
an interpretation in terms of the equivalence of two very different
looking categories.  The nature of the equivalence certainly
seems to equate odd D-branes (the Fukaya side) with even D-branes,
though the interpretation of D-branes as coherent sheaves (much
less a elements of the derived category) has as yet only been
speculative \hm, \mor.
 
We have offered a detailed proof of categorical mirror symmetry
in the case of the elliptic curve. However, our approach was
purely computational. A natural geometric construction of our
equivalence can be found in \ap.

One cannot at present be too optimistic about generalizing this
result to higher-dimensional Calabi-Yau manifolds.  The problem
is that the Fukaya category composition depends on the K\"ahler
form, i.e. on the unique Calabi-Yau metric.  The
exact form of this metric is of course
not known.  Abelian surfaces offer the only possible simple
extensions of this work -- this is work in progress.
Some hope may lie in the fact that areas of calibrated
submanifolds are topoligical, in a way, as the volume is the
restriction of a closed form.
 
Another obstacle is that the derived category of coherent
sheaves on a general Calabi-Yau manifold
is not easily studied.  For K3 some facts are known
\kthree, and through the work of \fmw\ we can expect some
simplification for elliptically fibered manifolds,
but a precise statement of the composition eludes
us.  One may hope to construct a map between objects,
however.  This may be a possible starting point.  For example,
the conjecture of \syz\ can be interpreted as simply positing an
object on the Fukaya side equivalent to the skyscraper sheaf
over a point on the mirror (derived category) side -- the
identification of the moduli space of this
object is immediately seen to be the mirror manifold itself.
Though Lagrangian submanifolds with flat $U(1)$ bundles
(the conjectured mirror objects to the skyscrapers) may not
seem too natural, we recall that the prequantum line bundle,
whose first Chern class equals the K\"ahler class, restricts
to a flat line bundle on Lagrangian submanifolds in the Calabi-Yau.
 
Even less is known about the Fukaya category.  Special-Lagrangian
submanifolds are not a well-studied class of objects.  Recent
work of Hitchin \hitch\ has only begun the investigation of
these spaces as dually equivalent to complex submanifolds.

The compositions in the Fukaya category contain
a wealth of information about holomorphic discs in the mirror
manifold, though it is not yet clear exactly how the many results of
ordinary mirror symmetry would follow from categorical mirror symmetry
for a general Calabi-Yau manifold. It seems that one has to start
with the algebraic problem of effectively computing invariants of
$A_{\infty}$-categories.

Finally, we anticipate a return to the connection of \helix, in which
the derived category of coherent sheaves on a Fano variety
was used to count the number of solitions of a supersymmetric
sigma model on that manifold.  In that context, the moduli space
was the space of topological field theories, or in other words 
the ``enlarged''  moduli space.  Armed with the interpretation of
D-branes as solitonic boundary states in a superconformal field
theory, we learn through Kontsevich's work that
the appearance of the derived category, while still somewhat
mysterious, is not likely to be accidental.
 
We look forward to further research on these matters.
\vskip.15in
{\bf Acknowledgements}
  
The research of  A.P. is supported in part by
NSF grant DMS-9700458, and that of E.Z. by grant
DE-F602-88ER-25065.

\listrefs
 
\end